\documentclass[12pt]{article}
\newcommand{\comment}[1]{}

\newcommand{\complexes}{{\bf C}}

\renewcommand{\Re}{\mathop{\rm Re}\nolimits}

\newcommand{\qed}{\vrule height 6pt depth 0pt width 3 pt}

\newcommand{\BibTeX}{{\rm B\kern-.05em{\sc i\kern-.025em b}\kern-.08em     
    T\kern-.1667em\lower.7ex\hbox{E}\kern-.125emX}}

\newcommand{\off}[1]{#1^{\mbox{\tiny off}}}
\newcommand{\diag}[1]{#1^{\mbox{\tiny d}}}

\newenvironment{proof}[1][Proof]{\begin{trivlist}\item[\hskip \labelsep
{\it #1. }]}{\hfill \qed \goodbreak \end{trivlist}}   

\newtheorem{theorem}{Theorem}

\newtheorem{corollary}[theorem]{Corollary}
\newtheorem{lemma}[theorem]{Lemma}

\begin{document}

\title{Estimates for the scattering map associated to a
two-dimensional first order system}
\author{Russell M. Brown\footnote{Author supported in part
by the National Science Foundation, Division of Mathematical
Sciences.} \\Department of Mathematics\\University of 
Kentucky\\Lexington, KY 40506-0027\\
\\
rbrown@pop.uky.edu }
\date{}

\maketitle

\begin{abstract}
We consider the scattering transform for the first order
system in the plane,
$$
\left(\begin{array}{cc}\partial_{\bar x}&  0\\
0 & \partial _x \end{array}\right) \psi - 
\left( \begin{array}{cc} 0 & q^1 \\q^2 & 0 \end{array}
\right)  \psi = 0. 
$$ We show that the scattering map is Lipschitz continuous on a neighborhood of
zero in $L^2$.
\end{abstract}

This paper gives an estimate for the scattering map associated to a
first-order system 
\begin{equation}\label{sys}
D \psi - Q \psi = 0
\end{equation}
in the plane. Here, $D$ and $Q$ are defined by 
$$
D=\left( \begin{array}{cc}  \partial_{\bar x} & 0 \\ 0 & \partial_{x}
\end{array}\right) \quad \mbox{and} \quad 
Q = \left( \begin{array}{cc}0 & q^1 \\ q^2 & 0 \end{array}\right)
$$
and  $\partial_{\bar x}$ and $\partial _{x}$ are the standard
derivatives with respect to  $x= x^1+ix^2$ and $ \bar x$. The entries
of the matrix $Q$, $q^1(x)$ and $q^2(x)$ are complex valued functions
on the complex plane. (We will consistently use superscripts to
indicate components.)
 The system (\ref{sys}) was
studied by Beals and Coifman \cite{BC:1985,BC:1988}, and a number of other
authors (see Fokas and Ablowitz \cite{FA:1984} for an earlier formal
treatment, Sung \cite{LS:1994a,LS:1994b,LS:1994c} for a detailed
rigorous treatment  and  the review articles
\cite{MR90f:35171,MR93h:35177} for additional references). 
  One beautiful feature of the
scattering theory as
developed by  Beals and Coifman is the existence of a non-linear
Plancherel identity relating the scattering data $S$,  which is  defined
below, to the potential $Q$ in (\ref{sys}),
\begin{equation}\label{nlp}
\int_{\complexes} |S|^2\,d\mu = \int_{\complexes} |Q|^2\,d\mu.
\end{equation}
Here, we use $d\mu$ to denote Lebesgue measure and $2\times
2$-matrices are normed by  $|F|^2= \sum_{j,k=1}^2 |F^{jk}|^2$.
The identity (\ref{nlp}) is valid for nice potentials $Q$ satisfying
one of the symmetries 
$Q^*= \pm Q$.  Since the map  $Q\rightarrow S$ is non-linear, the
identity $(\ref{nlp})$ does not immediately imply the continuity of
this map. The goal of this paper is to provide a proof of the
continuity in the $L^2$-metric of 
the map $ Q\rightarrow S$
for small potentials.  Our methods apply equally to the inverse map
$S\rightarrow Q$. We also compute explicitly the size of potentials
for which our method works and obtain results for potentials with
$\|Q\|_2 < \sqrt 2$. 

Before describing the scattering theory and the proofs of our
estimates, we indicate several reasons for interest in this
system. The scattering theory for this system can be used to transform
 the Davey-Stewartson II system to a linear evolution
equation. Thus,  our estimates can be used to assert that for small
initial data, the solution $u$ depends continuously on the  initial data  in
the $L^2$-norm. However, it  is not clear that this method  produces a
solution for all initial data in a neighborhood of zero in  $L^2$.  
The application to the Davey-Stewartson system seems to have motivated
the work of Beals 
and Coifman, Fokas and Ablowitz and Sung.

This system also appears in the work of the author and G. Uhlmann
\cite{BU:1996} who study the inverse conductivity problem in two
dimensions and prove that the coefficient is uniquely determined by
the Dirichlet to Neumann map for conductivities which have one
derivative in $L^p$, $p>2$.  This is the least restrictive regularity
assumption on the coefficient which is known to imply a uniqueness
theorem for the coefficient. The argument depends in a 
crucial way on the non-linear Plancherel identity (\ref{nlp}).  More
recently, Barcel\'o, Barcel\'o and Ruiz \cite{BBR:1999} have shown that
the coefficient depends continuously on the Dirichlet to Neumann map
under the {\em a priori } regularity hypothesis that the coefficient
has a gradient which is H\"older continuous.  One step of their
argument is the observation that the scattering map is continuous if
the potential is assumed to be in $C^\epsilon$ and compactly
supported.  The research reported here is perhaps a step towards
 relaxing this assumption.

We begin the formal development by reviewing the notation of Beals
and Coifman's paper \cite{BC:1988}. 
 We let $x= x^1 + ix^2$ and $ z=
z^1+i z^2$ denote variables in the complex plane.
A family of solutions of the free system, $D\psi_0 = 0$,  is given by
$$
\psi_0 (x,z) = \left(\begin{array}{cc} \exp ( ixz) & 0 \\ 0 &
\exp(-i\bar x z)     \end{array}\right).
$$
We look for solutions of (\ref{sys}) of the form $ \psi(x, z) =
 m(x,z)\psi_0(x,z)$ where $m$ approaches the $2\times 2$ identity
matrix as $|x|\rightarrow \infty$.  A computation shows that $m$
must satisfy the equation
\begin{equation}\label{sysz}
D_z m =Qm
\end{equation}
where $D_z$ is the operator $E_z^{-1} D E_z$ and for $z\in
\complexes$, $E_z$ is the map on $ 2\times 2 $ matrix-valued functions
defined by  
$$
E_z f = \diag{f} + A_z^{-1} \off{f}.
$$
Here and below, we use $\diag{f}$ and $\off{f}$
to denote the diagonal and off-diagonal parts of the matrix $f$ and
$A(x,z)= A_z(x)$ is the matrix given by
$$
A(x,z) = \left( \begin{array}{cc} \exp( i x \bar z + i \bar x z ) & 0
\\ 0 & \exp( -ixz -i\bar x\bar z ) \end{array}\right) 
= \left( \begin{array} {cc} a^1(x,z) & 0 \\ 0 & a^2(x,z)\end{array}\right).
$$
The solution to (\ref{sysz}) is found by solving the integral equation 
\begin{equation}\label{inteqn} 
m = 1 + G_z( Qm).
\end{equation}
In this equation and below, we use $1$ to denote the $2\times 2$
identity matrix. 
We use $G_z = E_z^{-1} G E_z$  to denote the inverse of $D_z$, where
$$
G(f)(x) = \frac 1 \pi\int_{\complexes } \left(\begin{array}{cc} (x-y)^{-1}
& 0 \\ 0 & ( \bar x - \bar y ) ^{-1} 	      \end{array} \right) f(y) \, d\mu(y) .
$$
 In our situation, with $Q$ small, we can write the solution as a
Neumann series, $m= \sum_{k=0} ^ \infty (G_zQ)^k (1) $.  Here and   throughout
this paper, we use $Q$ to denote both the function $Q$ and
multiplication  operator $ f \rightarrow Qf$.  

Continuing with the scattering theory, we can differentiate the
solution with respect to $z$ and obtain the $\partial_{\bar
z}$-equation
\begin{equation}\label{dbar}
\partial _{\bar z } m(x,z)= (Tm)(x,z)
\end{equation}
where the right-hand side $Tm$ is defined by 
\begin{equation} \label{tdef}
Tm(x,z) = m(x,\bar z) S(z) A(x,-\bar z).
\end{equation}
In the definition of $T$,  $S$ is the scattering data which is  defined by
\begin{equation} \label{scatdef}
S(z) = - \frac 1 \pi {\cal J} \int_{\complexes} E_z ( Q(x) m(x,z))\,
d\mu(x)
\end{equation}
and  ${\cal J}f = 2 \off{(Jf)}$ with $J$ the $2\times 2$ matrix,
$$
J = \frac 1 2 \left( \begin{array}{rr} -i & 0 \\ 0 & i \end{array}\right).
$$
The scattering data $S$ defined in (\ref{scatdef}) is the function
which appears in the left-hand side of (\ref{nlp}).
The identity (\ref{nlp}) is proven for nice functions  by Beals and Coifman in
\cite{BC:1988}. A detailed argument is given by Sung in 
\cite{LS:1994a,LS:1994b,LS:1994c} for potentials $Q$ which are in
$L^1\cap L^\infty$ and much of his work extends to potentials in
$L^p\cap L^{p'}$ ($p\neq 2$).  

Our main result is to show that for $Q$ small in $L^2$, (and say $Q\in
\cal S(\complexes)$, the Schwartz space), the map $Q\rightarrow S (Q)$
is Lipschitz continuous in $L^2$ and that a similar result holds for
the inverse map 
$S\rightarrow Q(S)$.  These results allow us to extend the maps to a
neighborhood of zero in $L^2$.  In fact, the scattering map
$Q\rightarrow S(Q)$ and the inverse map $S\rightarrow Q(S)$ are of the same
form. (The more adroit normalization used by Sung makes this clear.)
Thus, we concentrate most of our energy on the scattering map.

We remark that the definition of $S$ given in (\ref{scatdef}) does not
make sense if $Q \in L^2$, since I do not know if we can construct $m$
without the 
assumption of some extra decay on $Q$. One assumption that will
suffice is that $Q \in L^p$ for some $p<2$.  

Our proof is straightforward.  We write $m$, the solution of
(\ref{inteqn}) as a Neumann series and substitute this series into the
definition of $S$, (\ref{scatdef}).  This gives a series of
multilinear expressions in $Q$. We estimate each term and show that we
can sum the series.  The proof of these estimates depends only on the
Hardy-Littlewood-Sobolev theorem, but requires a certain amount of
persistence.    

The main ingredient is the following well-known result on the first
order Riesz potential defined by 
$$
R(f) (x) = \int_{\complexes} \frac {f(y)}{|x-y|} \, d\mu(y).
$$
We will only consider $R$ acting on scalar functions. 

\begin{lemma}\label{riesz} Let $r$ and $p$ satisfy $r>0$ and  $1<p/r<2$, then the 
map $f\rightarrow [R(|f|^r)]^{1/r}$ satisfies
$$
\|[ R(|f|^r)]^{1/r}\|_{\tilde p} \leq \alpha(r/p)^{1/r} \| f\| _p
$$
where $ 1/p - 1/\tilde p = 1/(2r)$.  The constant $\alpha(\theta)$ is
bounded  for  $\theta$  in compact subsets of the interval
$(1/2,1)$ and we have
$$
\alpha(3/4) = \pi \qquad \mbox{and}\qquad 
\limsup _{\theta\rightarrow 3/4}  \alpha (\theta) \leq \pi.
$$
\end{lemma}

\begin{proof} According to the Hardy-Littlewood-Sobolev theorem, on
fractional integration, (see \cite{ES:1970}, for example), we have
$R: L^s \rightarrow L^{ \tilde s}$ where $ 1/s-1/\tilde s= 1/2$,
provided $ 1< s< 2$. We let $ \alpha(1/s)$ denote the operator norm of
this map on scalar-valued functions.

Now, given $p$ and $r$  as in the Lemma, we define $ \tilde p$ by
$r/p - r/\tilde  p = 1/2$. We apply the
Hardy-Littlewood-Sobolev theorem with exponents $p/r$ and $ \tilde
p/r$  to conclude that 
$$
\left( \int_{\complexes}[R(|f|^r) (x)] ^{\tilde p/r }
\,d\mu(x)\right)^{r/{\tilde p}} \leq \alpha(r/p) \left( \int_{\complexes} 
|f(x)|^p \, d\mu(x)\right)^{r/p}.
$$
Taking the $r$th root gives the first inequality in the Lemma. 

The operator norm for $1/s= 3/4$ was computed by Lieb
\cite{MR86i:42010}.  
The behavior 
near $\theta = 3/4$ asserted in the Lemma follows from Lieb's
result, the boundedness for $ 1<p<2$  and the Riesz-Thorin
interpolation theorem.  
\end{proof}

We now state the main result.
\begin{theorem} \label{main}
Let 
$$
N= L^2 (\complexes) \cap \{ F : \|F\|_2 <\sqrt 2 \mbox{ and } F^d=0\}
$$
then the maps $S\rightarrow Q$ and $Q \rightarrow S$, defined initially
on ${\cal S} (\complexes)\cap N$, extend continuously to $N$.  
In addition, 
$$ S\circ Q = \mbox{Id}\quad \mbox{and} \quad Q\circ S = \mbox {Id}
$$
provided $Q$ and $S(Q)$ lie in $N$ (respectively, $S$ and
$Q(S)$ lie in $ N$). If $ Q = \pm Q^*$, then
$$
\int_{\complexes} |Q|^2\,d\mu = \int_{\complexes} |S|^2\,d\mu.
$$
\end{theorem}

The proof will depend on multi-linear estimates involving repeated
fractional integration.  For the study of these multi-linear
expressions, we will work in a sequence of $L^p$ spaces using the
exponents, $ p_j$, $s_j$ and $r_j$ defined by:
$$ 
p_0=2 \mbox { and } s_0 =4.
$$
Then we define
\begin{eqnarray} 
\frac 1 { r_j } & =& \frac 4 { 3 p _j } , \qquad j = 0,1,2, \dots
\label{FirstExponent} \\ 
\frac 1 {p_{ j+1} } & =& \frac 1 { p_j} - \frac 1 { 2r_j }  , \qquad
j =0, 1, 2,\dots \\
\frac 1 { s_{j+1} } & =& \frac 1 { s_j } +\frac 1 { 2r_j}, \qquad
j=0,1,2\dots \\ 
\frac 1 {\tilde s_{ j} }  & = & \frac  1 { s_{ j} } - \frac 1 2 
 \qquad j = 1, 2, \dots \label{LastExponent}
\end{eqnarray}
These sequences satisfy
\begin{eqnarray} 
\frac { p_j } { r_j } & = & \frac 4 3, \qquad j=0,1,2,\dots
\label{PoverR} \\
\frac 1 { p_j} + \frac 1 { s_j} &  =& \frac 3 4  , \qquad j = 0, 1,2,
 \dots \label{ConservationLaw} \\
 \frac { s_j } { r'_j } & = & \frac 4 3, \qquad j=
0,1,2,\dots \label{PprimeoverR}  \\
 s_{ j} & < &  2, \qquad j  =  1,2,\dots \label{Sinequality}
\end{eqnarray}
Here and below, we use $ r'= r/(r-1)$ to denote the dual exponent. 
The statement (\ref{PoverR}) is  immediate from
(\ref{FirstExponent}) and (\ref{ConservationLaw}) follows from the
definitions of $p_j$ and $s_j$.  The third,
(\ref{PprimeoverR}), follows since (\ref{ConservationLaw}) implies
 $ s_j = 4p_j/ ( 3p_j -4)$ and $ 1/r_j' = ( 3p_j -4)/ (3p_j) $. The last
observation (\ref{Sinequality})  follows since $s_1 = 12/7$ and $s_j$
decreases.

The multi-linear expression we will study is 
$I_k ( t, q_0, \dots q_{ 2k})$ which is defined  for $k =1,2,\dots$ by 
\begin{eqnarray*}
\lefteqn{
I_k( t, q_0,\dots, q_{2k})
} & & \\
& =&  \int_{\complexes ^{2k+1} }\frac { t 
(x_0-x_1+ x_2 -\dots -x_{2k-1}+ x_{2k}) \prod_{j = 0 } ^{2k} q_j (
x_j)}{ |x_0-x_1||x_1-x_2|\cdots |x_{2k-1}- x_{2k}| } \, d\mu(x_0,\dots
x_{2k}).
\end{eqnarray*}
Also, we will use $I_0( t,q_0) = \int _{\complexes} t(x) q_0(x) \,
d\mu(x)$. 
The functions $t$ and $q_j$ are scalar-valued. 
The main estimate for this expression  is in the following Lemma. 
\begin{lemma} \label{mainlemma}
For every $\epsilon>0$, there exists a constant $C_\epsilon$ so that 
$$
I_k(t, q_0,\dots, q_{2k}) \leq C_\epsilon \pi ^{2k} (1 + \epsilon)^{2k}
\|t\|_2 \prod_{j=0}^{2k} \| 
q_j\|_2.
$$
\end{lemma}

Before, beginning the proof of  Lemma \ref{mainlemma},  we
state and prove the following result. The proof of  Lemma
\ref{mainlemma} amounts to applying this result $k$ times. 

\begin{lemma} \label{inductivestep}   Let $j \geq 0$,  $k \geq 1$
 and suppose that  $p_j$,
$s_j$, $ \tilde s_{ j+1}$   and $r_j$ are exponents as defined in
(\ref{FirstExponent}-\ref{LastExponent}). 
For $ t $ and $ q_j$ non-negative functions on $ \complexes$, we have 
$$
I_ k ( t, q_0, \dots, q_{2k} )\leq I_{ k-1} ( t_1, q_2\tilde q_2, q_3
, \dots , q_{ 2k}) 
$$
where  the functions  $t_1$ and $\tilde q_2$ satisfy
\begin{eqnarray*}
\| t_1\| _{p_{j+1}} & \leq& \alpha (3/4) ^{ 1/r_j} \| t\|_{ p_j} \\
\| \tilde q_2\| _{ \tilde s_{j+1}} & \leq & \alpha (3/4 )^{ 1/ r'_j}
\alpha ( s_{j+1} ^ { -1})\| q_{1}\| _ 2 \| q_0 \|_{ s_j}.
\end{eqnarray*}
\end{lemma}

\begin{proof} We consider the integral with respect to $x_0$  in $I_k$
and apply the H\"older inequality using exponents $ r_j$ and $r'_j$
to obtain
\begin{eqnarray*}
\lefteqn {\int _{ \complexes}\frac{ t (x_0-x_1+ x_2-\dots + x_{
 2k})q_0 ( x_0) }{|x_0-x_1|} \, d\mu (x_0)}\\ & \leq & \left( \int _{
 \complexes}\frac{ q^{r'_j}_0 ( x_0) }{|x_0-x_1|} \, d\mu (x_0)
 \right)^{ 1/r'_j} \\
 & & \qquad \times \left( \int _{
 \complexes}\frac{ t^{r_j} (x_0-x_1+ x_2-\dots + x_{ 2k}) }{|x_0-x_1|}
 \, d\mu (x_0) \right)^{ 1/r_j} \\ & = &   \tilde q_1(x_1) \,  t_1(x_2- x_3+ x_4- \dots + x_{
 2k}).
\end{eqnarray*}
Here, the function $t_ 1$  is defined by 
$$
t_1( x) = [ R(t^{ r_j} )]^{ 1/r_j}(x)
        =  \left (  \int _{ \complexes }  \frac { t^{ r_j} ( x-w)}
        {|w|} d\mu(w) 
        \right)^{ 1/ r_j } .
$$
One may make the change of variables $ w= x_1-x_0$ in the integral
with respect to $x_0$ to obtain this representation of $t_1$. 
The function $ \tilde q_1$ is defined by 
$
\tilde q_1 (x_1) =  [R(q_0^{r'_j})]^{1/{ r'_j}} (x_1).
$
Now, we can rewrite the  integral with respect to $x_1$ in  $I_k$ as 
$$
 \int_{\complexes}  \frac  { q_1( x_1) \tilde q_1( x_1) } { | x_2-
x_1|}\, d\mu(x_1)  
= R ( q_1 \tilde q_1)(x_2) = R ( q_ 1 [R( q_ 0^{ r'_j})]^{ 1/{r'_j} })
(x_2)\equiv \tilde q_2(x_2)
$$
which defines $ \tilde q_2$. 
Inserting the definitions of $ t_1$ and $ \tilde q_2$ into $I_k$ gives
the inequality relating $ I_k$ and $I_{ k-1}$. 

We now establish the estimates for $t_1$ and $\tilde q_2$. 
The estimate 
$$
\| t_ 1\| _ {p_{j+1}} \leq \alpha (3/4)^{ 1/ r_j} \| t\|_{ p_j}
$$
follows from Lemma \ref{riesz} with $p = p_j$ and $ r= r_j$. The
hypothesis $r_j/p_j \in ( 1/2, 1)$ follows from (\ref{PoverR}). 
To estimate $ \tilde q_2$,  observe that $ 4/3< s_{j+1}<2$ from
(\ref{ConservationLaw}) and (\ref{Sinequality}) and thus we may apply  Lemma \ref{riesz} with $ p
= s_{j+1}$ and 
$r=1$ to obtain the first inequality below. Next, we use  H\"older's inequality and  Lemma \ref{riesz} with $ p
= s_j$ and $ r= r'_j$ to obtain the second and third inequalities:
\begin{eqnarray*}
\| R( q_1[R( q_0^{ r'_j})]^{ 1/ r'_j}) \| _ {\tilde s_{j+1}} 
& \leq &  \alpha ( s_{j+ 1} ^{ -1} ) \| q_ 1 [R( q_ 0 ^{ r_j '}  ) ]
^{ 1/r_j' } \| _ { s_{j+1}}  \\
 & \leq &  \alpha ( s_ {j+1} ^{ -1} ) \| q_{1}\| _2 \| [R( q_0 ^{ r_j '}  ) ]
^{ 1/r_j' } \| _ {\tilde s_{j+1}}  \\
& \leq & \alpha (s_{j+1}^{ -1} ) \alpha (3/4 )^{ 1/ r_j ' }
 \| q_ 1\| _2 \| q_0 \| _{s_j} .
\end{eqnarray*}
These are the estimates of the Lemma. 
\end{proof}

We are now ready to give the proof of Lemma \ref{mainlemma}

\begin{proof}[Proof of Lemma \ref{mainlemma}] It suffices to prove the
Lemma when all of the functions are positive. We first claim that
if $ s= 4$ or $4/3$, we have
\begin{equation}\label{firststep}
|I_k(t, q_0, \dots, q_{ 2k}) | \leq C _\epsilon \pi^{ 2k} ( 1+
 \epsilon ) ^{ 2k} \|t\|_2 \| q_ 0\|_s \|q_{ 2k}\|_{ s'} \prod _{
 j=1}^{ 2k-1} \| q_j \| _2.
\end{equation}
We use the sequence of exponents defined in
(\ref{FirstExponent}--\ref{LastExponent}). By Lemma 
\ref{inductivestep} applied $k$ times, we obtain
$$
I_k(t, q_0, \dots, q_{ 2k}) \leq  I _0( t_k, q_{ 2k} \tilde q_{
 2k})
$$
where the function $t_k$ satisfies
$$
\|t_k \|_{ p_k } \leq \alpha( 3/4) ^{ \sum _{ j=0}^{ k-1} \frac 1 {
r_j } }\| t\|_{ p_0}
$$
and the function $ \tilde q_{ 2k}$ satisfies
$$
\| \tilde q _{ 2k}\| _{ \tilde s _ k} \leq \alpha( 3/4) ^{ \sum _{
j=0}^{ k-1} \frac 1 {r'_j } } \left( \prod _{ j=1}^k \alpha(s_j^{ -1}) \right)
\prod_{ j=1} ^ { 2k-1} \| q_j \|_2 \| q_0\|_{ s_0}.
$$
We write out the integral defining $I_0$ and obtain
$$
I_0 ( t_k, q_{ 2k} \tilde q_{ 2k}) = \int_ {\complexes} t_k (x_{2k})
 q_{ 2k}( x_{2k}) \tilde q _{ 2k}( x_{ 2k})\, d\mu(x_{ 2k})
$$
We have $ t_k \in L^{ p_k}$, $q_{2k} \in L^{ s'_0}$ and $ \tilde
q_{2k} \in L^{ \tilde s_{k}}$ where  the exponents $p_k$, $\tilde s_k
$ and $ s'_0$ satisfy  
$$
\frac 1 { p_k } + \frac 1 { \tilde s _k } + \frac 1 { s_0 '} = 
\frac 1 { p_0 }  + \frac 1 { s_0 } -\frac 1 2 + \frac 1 { s_0'}= 1.
$$
Thus H\"older's inequality and the estimates for $ t_k$ and $ \tilde
q_{2k}$ imply that  
$$
| I_k ( t, q_0, \dots, q_{2k} ) | \leq \alpha(3/4)^k \prod_{ j=1}^k
 \alpha( s_j ^{ -1})  \| q_0 \|_4 \| q_{ 2k}
  \|_{ 4/3}\|t\|_2 \prod_{ j=1}^{ 2k-1} \| q_j \| _2 .
$$
Since $s_j \rightarrow 4/3$ as $j \rightarrow \infty$,  Lemma
\ref{riesz}  tells us that $\limsup _{ 
j\rightarrow \infty }  \alpha(s_j ^{ -1}) \leq \pi$. Thus,  if $
\epsilon > 0$, then we obtain (\ref{firststep}) with $ s=4$.  To
obtain the result with $ s=4/3$, we simply need to change 
variables to replace  $ x_j $ by $ x_{ 2k-j}$. 

Now that we have the estimate (\ref{firststep}), 
we apply a simple interpolation argument to obtain the
Lemma. We define an operator ${\cal T}$ by 
$$
{\cal T} q_{ 2k}( x_0) = \int _{ \complexes ^{ 2k}}  \frac { t( x_0-x_1 +
\dots + x_{ 2k}) q_ 1( x_1) \dots q_{ 2k } (x_{2k})} { | x_0-x_1|\dots
|x_{2k-1} - x_{2k}|}\, d\mu(x_1,\dots, x_{ 2k}).
$$
According to the estimate (\ref{firststep}),  $ \cal T$ maps $L^4
\rightarrow L^4$ and $L^ {4/3} \rightarrow L^{ 4/3}$.  Hence, by the
Riesz-Thorin interpolation theorem, ${\cal T}$ is bounded on $L^2$, with the
same bound. This implies the Lemma. 
\end{proof}

\begin{proof}[Proof of Theorem 2]
We begin by indicating how we construct the solutions $m$ when the potential
$Q$ is nice, $Q\in {\cal S (\complexes)}$, say, and
$\|Q\|_2< \sqrt 2$. We can write $m$ as the infinite sum 
\begin{equation}\label{jost}
m = 1 +  \sum _ {j=1} ^\infty (G_zQ)^j (1).
\end{equation}
We consider one entry in the matrix $(G_zQ)^j(1)$ with  $j= 2k$. Note
that for $j$ even, the off-diagonal part $ (G_zQ)^j(1)$ is zero. Also,
note that the operator $E_z$ involves multiplying by exponentials of
modulus 1, and thus the entries of $E_zf$ have  the same $L^p$-norm as
the entries of $f$.  We apply Lemma \ref{riesz} with $p =4/3$ and $
r=1$ to obtain that  
\begin{eqnarray*} \| ((G_zQ)^{2k})(1)^{11} \|_4
 & \leq& \pi^{-2k} \|R(Q^{12}(\dots R(Q^{21})\dots)) \|_4 \\
&\leq & ( \|Q^{12}\|_2\|Q^{21}\|_2) ^{k-1}\|Q^{12}\|_2\|Q^{21}\|_{4/3} \\
&\leq & \frac 1   {2 ^{k-1}} \|Q\|_2^{2k-1} \|Q\|_{4/3}.
\end{eqnarray*}
The last inequality uses the elementary inequality that
$\|Q^{12}\|_2\|Q^{21}\|_2 \leq \frac 1 2 ( \|Q^{12}\|_2^2 +
\|Q^{12}\|_2^2)$.  A similar argument holds for the remaining entries.
Thus the infinite sum in (\ref{jost}) will converge in $L^4$ if
$\|Q\|_2 < \sqrt 2$.  We substitute the sum for $m$ into
(\ref{scatdef}). Since $Q$ and $S$ are off-diagonal, only the even
terms of the series for $m$ are needed in the expression for
$S$. Also, note that the operator $G_z= G$ when acting on diagonal
matrices. Using these observations and that $Q$ is off-diagonal, we
can write
\begin{eqnarray*}
S(z) & =& \frac {-2} \pi J \int_{\complexes} A(x,-z) Q(x)
\sum_{j=0}^\infty (GQG_zQ)^j(1)\, d\mu(x)\\
& \equiv & \frac {-2} \pi J \sum _{j=0}^\infty S_j(z).
\end{eqnarray*}

The term with $j=0$  in this series is essentially the Fourier transform. If
we consider one entry, we obtain
$$
S^{12}_0(z) = \int_{\complexes} \exp (- 2i ( x^1z^1+ x^2 z^2)) Q^{12}
(x) \, d\mu (x)
$$
and we have a similar expression for $S^{21}_0$. Thus,
$$\frac 1 {\pi ^2} \int_{\complexes} |S_0|^2\, d\mu = 
\int_{\complexes} |Q|^2 \, d\mu 
$$
by the standard Plancherel identity. The higher order terms require a
bit more work.  

In order to simplify the notation, we consider the upper-right entry
in $S_k$
\begin{eqnarray*}
S_k ^{12} (z) & =& \frac  {1 } {\pi ^{2k}}\int_{\complexes ^{2k+1}}
a^1(-x_0+ x_1-x_2+x_3\dots+x_{2k-1} - x_{2k},z)\\
 & & \quad \times  \frac { Q^{12}(x_0)
Q^{21}(x_1)\dots Q^{21}(x_{2k-1} )Q^{12}(x_{2k})}{ ( \bar x_0-\bar x_1)
( x_1-x_2) \dots (\bar x_{2k-2}- \bar x_{2k-1})(
x_{2k-1}-x_{2k})}\,d\mu(x_0,\dots , x_{2k}).
\end{eqnarray*}
We take a sufficiently regular (scalar-valued) test function $T$ and write out the
expression for $S^{12}_k$ to obtain
\begin{eqnarray*}
\int_{\complexes} T(z) S_k ^{12}(z) \, d\mu (z)& = &
 \frac  {1  } {\pi ^{2k}}\int_{\complexes ^{2k+2}}
T(z) a^1(-x_0+ x_1-x_2+x_3\dots+x_{2k-1} - x_{2k}, z) \\
\lefteqn{\hspace {-0.5 in}
\times  \frac { Q^{12}(x_0)
Q^{21}(x_1)\dots Q^{21}(x_{2k-1} )Q^{12}(x_{2k})}{ ( \bar x_0-\bar x_1)
( x_1-x_2) \dots (\bar x_{2k-2}- \bar x_{2k-1})(
x_{2k-1}-x_{2k})}\,d\mu(z, x_0,\dots , x_{2k}). } & &
\end{eqnarray*}
Note that if each of the  functions $T$ and $Q$ are in $L^1\cap
L^\infty$, say, then the integral on the right of the previous
expression converges absolutely. Thus, we may use the Fubini theorem
to carry out the integration with respect to $z$ first and obtain
the expression
\begin{eqnarray*}
\lefteqn{ \frac  {1 } {\pi ^{2k}}\int_{\complexes ^{2k+1}}
\hat T(2(x_0- x_1+ x_2 -\dots - x_{2k-1} + x_{2k})) } & & \\
 & &  \quad \times  \frac { Q^{12}(x_0)
Q^{21}(x_1)\dots Q^{21}(x_{2k-1} )Q^{12}(x_{2k})}{ ( \bar x_0-\bar x_1)
( x_1-x_2) \dots (\bar x_{2k-2}- \bar x_{2k-1})(
x_{2k-1}-x_{2k})}\,d\mu(x_0,\dots , x_{2k}).
\end{eqnarray*}
Here, $ \hat T$ is the Fourier transform of $T$ and is defined by
$\hat T ( z) = 
\int f(x) e^ {- \Re x\bar z}\, d\mu(x) = \int f(x) a^1( -\frac 1 2 x,
z)\, d\mu(z)$.  Thus we have that
$$
\left | \int_\complexes T(z) S^{ 12}(z) \, d\mu(z) \right |
\leq \frac 1 { \pi ^{2k}} I_k( | \hat T(2\cdot)|, | Q^{ 12}|,
|Q^{21}|,\dots , |Q^{12}|).
$$
Now, Lemma \ref{mainlemma} implies for each $\epsilon > 0$, there is a
constant $C_\epsilon$ so that 
$$
\|S_k\|_2 \leq  C_{\epsilon} \frac { (1+\epsilon)^{2k}} {2^k} \|Q\|^{2k}_2
$$ 
and hence  the series 
$$ \sum_{ k=0}^\infty S_k
$$ 
converges in $L^2$ if
$\|Q\|_2 < \sqrt 2$.  

We now turn to the proof that the map $ Q \rightarrow S$ is
Lipschitz continuous. Let $Q$ and $\tilde Q$ be two matrix potentials. We write
$S= S(Q)$, then we have
\begin{eqnarray*}
\lefteqn{\left| \int_{\complexes} T(z)( S_k ^{12}(Q)(z) - S_k^{12}(\tilde
Q)(z))\, d\mu(z)\right| } \\
&   \leq &\frac 1 {\pi^{2k}}  \sum_{j=0} ^{2k} I_k(| \hat
T(2\cdot)|, |Q^{12}|, 
|Q^{21}|,\dots, |Q^{ab}-\tilde Q ^{ab}|,\dots, | \tilde Q^{12}|)
\end{eqnarray*}
where the difference occurs in the $j$th spot and $ab= 12 $ if $j$ is
even or $21$ if $j$ is odd. 
Thus, Lemma \ref{mainlemma}  implies that if $\|Q\|_2<M$ and  $ \| \tilde Q
\|_2<M$, then  for each $\epsilon>0$, there is  a constant
$C_\epsilon$ so that 
$$
\|S^{12}(Q)-S^{12}(\tilde Q)\| _2 \leq C_\epsilon  \|Q-\tilde Q\|_2 M^2
\sum_{k=0}^\infty (2k+1)(1+\epsilon)^{2k}\frac 1 {2^{k-1}} M^{2(k-1)} .
$$
Since we have $M<\sqrt 2$, we have  the Lipschitz  continuity also.

The inverse map is handled similarly so we only give a sketch of the
argument.  The map which takes $Q$ to $S$  is given by 
\begin{equation}\label{s2q}
Q ( x) = \frac 1 \pi  {\cal J} \int _ {\complexes } Tm(x,z)\, d\mu(z) 
\end{equation}    
where the map $T$ depends on $S$ and was defined in (\ref{tdef}).
We have that $m$ may be represented as
$$
m = ( I-CT)^{-1} ( 1)
$$ 
where $C$ is the Cauchy transform
$$
Cf(z) = \frac 1 \pi \int _{\complexes} \frac { f ( w)}  { z-w}  \,
d\mu (w) 
$$
acting on matrix valued functions.  
As before, we can write $m$ as the series,
\begin{equation}\label{jost2}
m = \sum_{ j=0}^\infty (CT)^j (1)
\end{equation}
Substituting (\ref{jost2})
into (\ref{s2q}) gives a series representation for $Q$, $Q(x) =
\frac 1 \pi \sum_{k=0} ^\infty Q_k(x)$ where the $k$th term is given by 
\begin{eqnarray*}
Q_k ^{12} (x) &= & 
\frac {1 } {\pi ^{2k} }\int_{\complexes ^{2k+1}}
\frac { S^{12}( z_0) S^{21} ( z_1)\dots S^{12} ( z_{2k}) }
{ ( \bar z_0-z_1)(\bar z_1 - z_ 2) \dots ( \bar z _{2k-1}- z_{2k})} \\
 & & \qquad  \qquad \times  a^1( x, z_ 0- \bar z _1 + \dots - \bar z
_{ 2k- 1} + z_ {2k}) \,  d\mu ( z_ 0, \dots z_{2k}).
\end{eqnarray*}
Arguing as before, we can show that  for each $ \epsilon >0$, there is
constant $C_\epsilon$ so that 
$$
\| Q_{2k}  \|_2 \leq  C_\epsilon\frac { (1+\epsilon)^{2k}}{2^k}  \|
S\|_2 ^{2k+1} 
$$
and also obtain that $Q$ depends continuously on $S$.  

The statements that $Q\circ S = Id$ and $S\circ Q= Id$ follow from
continuity and the corresponding results for nice potentials as proved
by Sung \cite{LS:1994a,LS:1994b,LS:1994c}.  However, one must replace 
Lemma 2.3 of \cite{LS:1994a}. The uniqueness of  $m$ 
needed to carry Sung's arguments
 holds for  $ \|Q\|_2 < \sqrt 2 $, thanks to
the estimate for $ \alpha(3/4)$ in Lemma \ref{riesz} (see also
(\ref{jost})). Also, the extension of the 
non-linear Plancherel identity to functions in $L^2$ is immediate from
the result for smooth functions and the continuity of the map $
Q\rightarrow S$.  
\end{proof}

We give an application of our estimate. This application depends
on the fact that we have shown that the scattering map and its inverse 
map  a neighborhood of 0 in $L^2$ into $L^2$. Since the domain and
range lie in the same space, $L^2$, it is particularly simple to
construct the solution. 

Recall that the Davey-Stewartson II equation as studied by Beals and
Coifman is the evolution equation in two space dimensions given by
$$
\left\{ \begin{array}{l}
q_t = i q_{x^1x^2} - 4irq \\
r_{x^1x^1} + r_{x^2x^2} =  (|q|^2)_{x^1x^2}
\end{array} \right.
$$
If we set 
$$
Q (x,t) =\left( \begin{array}{rr}  0 &q(x,t) \\ \bar q(x,t) & 0 \end{array}\right)
$$
and apply the scattering map to obtain $S((\cdot, t)= S(Q)(\cdot, t)$, then the
function $S(t)$
satisfies
$$
S_t= (4 i z^1z^2) S. 
$$
Thus the solution of the Davey-Stewartson II equation is 
the upper-right entry of the
matrix
\begin{equation} \label{dssoln}
Q(\cdot, t) = Q( \exp(4itz^1z^2)S(Q(\cdot,0))).
\end{equation}
This is established for nice potentials $Q$ by Beals and
Coifman \cite{BC:1988}, see also Sung for a detailed argument
\cite{LS:1994a,LS:1994b,LS:1994c}.  Our contribution is the
observation that the 
maps $Q$ and $S$ are continuous in $L^2$ which immediately implies the
following result.   See Ghidaghlia and Saut \cite{GS:1990} and Linares
and Ponce \cite{LP:1993}  for additional results on the
Davey-Stewartson systems.

\begin{corollary} The scattering solution of the Davey-Stewartson II
system as defined in (\ref{dssoln}) gives a solution which depends
continuously on the initial data in $L^2$. If the initial data $q_j(\cdot,0)$ 
satisfies $\|q_j(\cdot,0)\|_2 <1$  
and say $q_j(x,0)\in {\cal S}(\complexes)$ for $j=1$ and 2, then  the
solutions $q_j(x, t)$ 
satisfy 
$$\| q_1(\cdot,t)- q_2(\cdot,t) \| _2 \leq C\| q_1 (\cdot ,0) -
q_2(\cdot ,0)\| _2.
$$
\end{corollary}

We close with three  questions related to the above results.
\begin{enumerate}
\item  Can we remove the smallness restriction from our Theorem
\ref{main}? It is known that if one of  $ Q = \pm Q^*$ holds  and $Q\in L^2$,
then the operator $I-G_zQ$ is invertible.  This fact does not seem to
be in the literature. To see the invertibility of $ I-G_z Q$ on $L^{4/3}$, we
observe that 
arguments similar to those of Nachman  \cite[Lemma 4.2]{AN:1994} show
that the map $ f \rightarrow G_z Q f$ is compact on $L^{4/3}$. The
injectivity of the map $ f \rightarrow f- G_zQf$ when $ Q= \pm Q^*$  follows from
Corollary 3.11 of \cite{BU:1996}. This Corollary is an observation of
Nachman's. Hence, the Fredholm theory gives the claim.

\item Can we construct the Jost solutions $m$ when $Q$ is only assumed
to be in $L^2$?

\item The expression (\ref{dssoln})   is defined,  if   $Q$ is just in
$L^2$ and small. In what sense does this expression solve the
Davey-Stewartson equations? 
\end{enumerate}


\begin{thebibliography}{10}

\bibitem{BBR:1999}
Juan~Antonio Barcel\'o, Tomeu Barcel\'o, and Alberto Ruiz, \emph{Stability of
  the inverse conductivity problem in the plane}, Preprint, 1999.

\bibitem{BC:1985}
R.~Beals and R.R. Coifman, \emph{Multidimensional inverse scatterings and
  nonlinear partial differential equations}, Pseudodifferential operators and
  applications (F.~Treves, ed.), Proceedings of symposia in pure mathematics,
  vol.~43, Amer. Math. Soc., 1985, pp.~45--70.

\bibitem{BC:1988}
\bysame, \emph{The spectral problem for the {D}avey-{S}tewartson and {I}shimori
  hierarchies}, Nonlinear evolution equations: Integrability and spectral
  methods, Manchester University Press, 1988, pp.~15--23.

\bibitem{MR90f:35171}
Richard Beals and R.~R. Coifman, \emph{Linear spectral problems, nonlinear
  equations and the $\overline\partial$-method}, Inverse Problems \textbf{5}
  (1989), no.~2, 87--130.

\bibitem{BU:1996}
Russell~M. Brown and Gunther~A. Uhlmann, \emph{Uniqueness in the inverse
  conductivity problem for nonsmooth conductivities in two dimensions}, Comm.
  Partial Differential Equations \textbf{22} (1997), no.~5-6, 1009--1027.

\bibitem{MR93h:35177}
A.~S. Fokas and L.-Y. Sung, \emph{On the solvability of the ${N}$-wave,
  {D}avey-{S}tewartson and {K}adomtsev-{P}etviashvili equations}, Inverse
  Problems \textbf{8} (1992), no.~5, 673--708.

\bibitem{FA:1984}
A.S. Fokas and M.J. Ablowitz, \emph{On the inverse scattering transform of
  multidimensional nonlinear equations related to first-order systems in the
  plane}, J. Math. Phys. \textbf{25} (1984), 2494--2505.

\bibitem{GS:1990}
J.M. Ghidaglia and J.C. Saut, \emph{On the initial value problem for the
  {D}avey-{S}tewartson system}, Nonlinearity \textbf{3} (1990), 475--506.

\bibitem{MR86i:42010}
Elliott~H. Lieb, \emph{Sharp constants in the {H}ardy-{L}ittlewood-{S}obolev
  and related inequalities}, Ann. of Math. (2) \textbf{118} (1983), no.~2,
  349--374.

\bibitem{LP:1993}
Felipe Linares and Gustavo Ponce, \emph{On the {D}avey-{S}ewartson systems},
  Ann. Inst. Henri Poincar{\'e} \textbf{10} (1993), 523--548.

\bibitem{AN:1994}
Adrian~I. Nachman, \emph{Global uniqueness for a two-dimensional inverse
  boundary value problem}, Ann. of Math. (2) \textbf{143} (1996), 71--96.

\bibitem{ES:1970}
E.M. Stein, \emph{Singular integrals and differentiability properties of
  functions}, Princeton University Press, Princeton NJ, 1970.

\bibitem{LS:1994a}
Li{-}Yeng Sung, \emph{An inverse scattering transform for the
  {D}avey-{S}tewartson {II} equations, {I}}, J. Math. Anal. Appl. \textbf{183}
  (1994), 121--154.

\bibitem{LS:1994b}
\bysame, \emph{An inverse scattering transform for the {D}avey-{S}tewartson
  {II} equations, {II}}, J. Math. Anal. Appl. \textbf{183} (1994), 289--325.

\bibitem{LS:1994c}
\bysame, \emph{An inverse scattering transform for the {D}avey-{S}tewartson
  {II} equations, {III}}, J. Math. Anal. Appl. \textbf{183} (1994), 477--494.

\end{thebibliography}

\providecommand{\bysame}{\leavevmode\hbox to3em{\hrulefill}\thinspace}

\smallskip \noindent AMS subject classification: 37K10

\smallskip \noindent Keywords:  Scattering theory, integrable
systems, Davey-Stewartson system

{\noindent \small \today}

\end{document}